\documentclass{amsart}
\usepackage{amssymb, amsmath, mathrsfs, amsthm}

\theoremstyle{plain}
\newtheorem{theorem}{Theorem}[section]
\newtheorem{lemma}[theorem]{Lemma}

\newtheorem{corollary}[theorem]{Corollary}

\newcommand{\MainTheoremName}{Main Theorem}

\theoremstyle{definition}
\newtheorem{definition}[theorem]{Definition}
\newtheorem{example}[theorem]{Example}

\theoremstyle{remark}
\newtheorem{remark}[theorem]{Remark}

\newtheorem{question}[theorem]{Question}

\numberwithin{equation}{section}

\newcommand{\nbd}{\nobreakdash}

\newcommand{\homeo}{\cong}

\newcommand{\powset}[1]{\mathcal{P}(#1)}
\newcommand{\card}[1]{\lvert #1\rvert}
\newcommand{\la}{\langle}
\newcommand{\ra}{\rangle}
\newcommand{\bd}{\partial}

\newcommand{\cardc}{\mathfrak{c}}
\newcommand{\cardp}{\mathfrak{p}}
\newcommand{\inv}[1]{#1^{-1}}
\newcommand{\closure}[1]{\overline{#1}}
\newcommand{\restrict}{\upharpoonright}
\newcommand{\subbase}[1]{\mathbb{S}(#1)}

\newcommand{\cell}[1]{c(#1)}
\newcommand{\weight}[1]{w(#1)}
\newcommand{\piweight}[1]{\pi(#1)}
\newcommand{\character}[1]{\chi(#1)}
\newcommand{\caliber}[1]{\mathrm{cal}(#1)}
\newcommand{\Aut}[1]{\mathrm{Aut}(#1)}

\newcommand{\reals}{\mathbb{R}}
\newcommand{\rationals}{\mathbb{Q}}
\newcommand{\Fn}[2]{\mathrm{Fn}(#1,\,#2)}

\newcommand{\mcA}{\mathcal{A}}
\newcommand{\mcU}{\mathcal{U}}
\newcommand{\mcV}{\mathcal{V}}
\newcommand{\msA}{\mathscr{A}}
\newcommand{\msB}{\mathscr{B}}
\newcommand{\msE}{\mathscr{E}}
\newcommand{\msR}{\mathscr{R}}
\newcommand{\msS}{\mathscr{S}}
\newcommand{\msT}{\mathscr{T}}

\DeclareMathOperator{\dom}{dom}

\DeclareMathOperator{\cf}{cf}
\DeclareMathOperator{\ind}{ind}

\begin{document}

\title{Amalgams, connectifications, and homogeneous compacta}
\author{David Milovich}
\address{University of Wisconsin-Madison Mathematics Department, 480 Lincoln Dr., Madison, WI 53706}
\email{milovich@math.wisc.edu}
\thanks{Support provided by an NSF graduate fellowship.}

\begin{abstract}
We construct a path\nbd-connected homogenous compactum with cellularity $\cardc$ that is not homeomorphic to any product of dyadic compacta and first countable compacta.  We also prove some closure properties for classes of spaces defined by various connectifiability conditions.  One application is that every infinite product of infinite topological sums of $T_i$ spaces has a $T_i$ pathwise connectification, where $i\in\{1,2,3,3\frac{1}{2}\}$.
\end{abstract}

\maketitle

\section{Introduction}
In \cite{Maurice}, M. A. Maurice constructed a family of homogeneous compact ordered spaces with cellularity $\cardc$.  All these spaces are zero\nbd-dimensional.  The cone over any of these spaces is path\nbd-connected but not homogeneous or ordered.  Indeed, it is easy to see that no compact ordered space with uncountable cellularity can be path\nbd-connected.  However, there is a path\nbd-connected homogeneous compactum with cellularity $\cardc$ which, though not an ordered space, has small inductive dimension $1$; we construct such a space by amalgamating copies of powers of one of Maurice's spaces together.  Moreover, this space is not homeomorphic to a product of dyadic compacta and first countable compacta.  To the best of the author's knowledge, there is only one other example \cite{vanMill} of a homogeneous compactum not homeomorphic to such a product, and that example's existence is independent of ZFC.

This amalgamation technique also can be used to construct new connectifications, where a connected (path\nbd-connected) space $Y$ is a \emph{connectification (pathwise connectification)} of a space $X$ if $X$ can be densely embedded in $Y$, and the connectification is \emph{proper} if the embedding can be chosen not to be surjective.  Whether a space has a connectification is uninteresting unless we restrict to connectifications that are at least $T_2$.   For a broad survey of connectification results, see \cite{wilson}.  Our focus will be on which a $T_2$ ($T_3$, $T_{3\frac{1}{2}}$, metric) spaces have $T_2$ ($T_3$, $T_{3\frac{1}{2}}$, metric) connectifications or pathwise connectifications.  Only partial characterizations are known.  For example, Watson and Wilson \cite{watson} showed that a countable $T_2$ space has a $T_2$ connectification iff it has no isolated points.  Emeryk and Kulpa \cite{emeryk} proved that the Sorgenfrey line has a $T_2$ connectification, but no $T_3$ connectification.  Alas {\it et al} \cite{alasetal} showed that every separable metric space without nonempty open compact subsets has a metric connectification.  Gruenhage, Kulesza, and Le Donne \cite{gruenhage} showed that every nowhere locally compact metric space has a metric connectification.

There are only a handful of results about pathwise connectifications.  For example, Fedeli and Le Donne \cite{fedelipath} showed that a nonsingleton countable first countable $T_2$ space has a $T_2$ pathwise connectication iff it has no isolated points.   Druzhinina and Wilson \cite{druzhinina} showed that a metric space has a metric pathwise connectification if its path components are open and not locally compact; similarly, a first countable $T_2$ ($T_3$) space has a $T_2$ ($T_3$) connectification if its path components are open and not locally feebly compact.  See also \cite{costantini} for some results about pathwise connectifications of spaces adjoined with a free open filter.

Suppose $i\in\{1,2,3,3\frac{1}{2}\}$ and $X$ has a proper $T_i$ connectification.  Then $X\times Z$ has a proper $T_i$ connectification for all $T_i$ spaces $Z$.  Thus, given one proper connectification, this product closure property gives us a new connectification.  We omit the easy proof of this fact here because we shall prove much stronger amalgam closure properties, which in many cases are also valid for pathwise connectifications.  The reals are a pathwise connectification of the Baire space $\omega^\omega$ because $\omega^\omega\homeo\reals\setminus\rationals$.  By applying amalgam closure properties to this particular connectification, we shall prove the following theorem.

\begin{theorem}\label{THMinfprodinfunion}
If $i\in\{1,2,3,3\frac{1}{2}\}$, then every infinite product of infinite topological sums of $T_i$ spaces has a $T_i$ pathwise connectification.  Every countably infinite product of infinite topological sums of metrizable spaces has a metrizable pathwise connectification.
\end{theorem}

The previously known result most similar to Theorem~\ref{THMinfprodinfunion} is due to Fedeli and Le Donne \cite{fedelioc}: a product of $T_2$ spaces with open components has a $T_2$ connectification iff it does not contain a nonempty proper open subset that is $H$\nbd-closed.

\section{Amalgams}

For all undefined notions, see \cite{engelking,juhasz}.

\begin{definition}
Given a topological space $X$, let $\subbase{X}$ denote the set of all subbases of $X$ that do not include $\emptyset$.
\end{definition}

Let $X$ be a nonempty $T_0$ space and let $\msS\in\subbase{X}$.  For each $S\in\msS$, let $Y_S$ be a nonempty topological space.  The \emph{amalgam} of $\la Y_S: {S\in\msS}\ra$ is the set $Y$ defined by
\begin{equation*}
Y=\bigcup_{p\in X}\prod_{p\in S\in\msS} Y_S.
\end{equation*}
We say that $X$ is the \emph{base space} of $Y$.  For each $S\in\msS$, we say that $Y_S$ is a \emph{factor} of $Y$.  Every amalgam has a natural projection $\pi$ to its base space: because $X$ is $T_0$, we may define $\pi\colon Y\rightarrow X$ by $\inv{\pi}\{p\}=\prod_{p\in S\in\msS} Y_S$ for all $p\in X$.  Amalgams also have natural partial projections to their factors: for each $S\in\msS$, define $\pi_S\colon\inv{\pi}S\rightarrow Y_S$ by $y\mapsto y(S)$.

Consider sets of the form  $\inv{\pi_S}U$ where $S\in\msS$ and $U$ open in $Y_S$.  We say such sets are \emph{subbasic} and finite intersections of such sets are \emph{basic}.  We topologize $Y$ by declaring these basic sets to be a base of open sets.  Let us list some easy consequences of this topologization.
\begin{enumerate}\label{amalgambasics}
\item For all $S\in\msS$, the map $\pi_S$ is continuous and open and has open domain.  

\item The map $\pi$ is continuous and open.

\item If $\card{Y_S}=1$ for all $S\in\msS$, then $Y\homeo X$.

\item For each $p\in X$, the product topology of $\prod_{p\in S\in\msS} Y_S$ is also the subspace topology inherited from $Y$.

\item Suppose, for each $S\in\msS$, that $Z_S$ is a subspace of $Y_S$.  Then the topology of the amalgam of $\la Z_S: S\in\msS\ra$ is also the subspace topology inherited from $Y$.

\item  Suppose, for each $S\in\msS$, that $\msS_S$ is a subbase of $Y_S$.  Then the set
\begin{equation*}
\{\inv{\pi_S}T: S\in\msS\text{\ and\ }T\in\msS_S\}
\end{equation*}
is a subbase of $Y$.
\end{enumerate}

Up to homeomorphism, an amalgam is a quotient of the product of its base space and its factors.  Specifically, the map from $X\times\prod_{S\in\msS}Y_S$ to $Y$ given by
\begin{equation*}
\la x,y\ra\mapsto y\restrict\{S\in\msS: x\in S\}
\end{equation*}
is easily verified to be a quotient map.

We say that a class $\mcA$ of nonempty $T_0$ spaces is \emph{amalgamative} if an amalgam is always in $\mcA$ if its base space and all its factors are in $\mcA$.  Therefore, any class of nonempty $T_0$ spaces closed with respect to products and quotients is amalgamative.  In particular, amalgams preserve compactness, connectedness, and path\nbd-connectedness. The next theorem says that several other well\nbd-known productive classes are also amalgamative.

\begin{theorem}\label{THMamalgamative}
The classes listed below are amalgamative provided we exclude the empty space.  Conversely, if an amalgam is in one of these classes, then its base space and all its factors are also in that class.
\begin{enumerate}
\item\label{T0amalgamative} $T_0$ spaces

\item\label{T1amalgamative} $T_1$ spaces

\item\label{T2amalgamative} $T_2$ spaces

\item\label{T3amalgamative} $T_3$ spaces

\item\label{T35amalgamative} $T_{3\frac{1}{2}}$ spaces

\item\label{hereddisconnamalgamative} hereditarily disconnected $T_0$ spaces

\item\label{zerodimamalgamative} zero\nbd-dimensional $T_0$ spaces
\end{enumerate}
\end{theorem}
\begin{proof}
For (\ref{T0amalgamative})\nbd-(\ref{T2amalgamative}), suppose $y_0$ and $y_1$ are distinct elements of $Y$.  If $\pi(y_0)=\pi(y_1)$, then there exists $S\in\dom y_0=\dom y_1$ such that $y_0(S)\not=y_1(S)$; whence, if $U_0$ and $U_1$ are neighborhoods of $y_0(S)$ and $y_1(S)$ witnessing the relevant separation axiom for $y_0(S)$ and $y_1(S)$, then $\inv{\pi_S}U_0$ and $\inv{\pi_S}U_1$ witness the the same separation axiom for $y_0$ and $y_1$.
If $\pi(y_0)\not=\pi(y_1)$, then let $U_0$ and $U_1$ be neighborhoods of $\pi(y_0)$ and $\pi(y_1)$ witnessing the relevant separation axiom for $\pi(y_0)$ and $\pi(y_1)$.  Then $\inv{\pi}U_0$ and $\inv{\pi}U_1$ witness the same separation axiom for $y_0$ and $y_1$.

For (\ref{T3amalgamative}) and (\ref{T35amalgamative}), suppose $C$ is a closed subset of $Y$ and $y\in Y\setminus C$.  Then there exist $n<\omega$ and $\la S_i\ra_{i<n}\in(\dom y)^n$ and $\la U_i\ra_{i<n}$ such that $U_i$ is a neighborhood of $y(S_i)$ for all $i<n$ and $\bigcap_{i<n}\inv{\pi_{S_i}}U_i$ is disjoint from $C$.   For each $i<n$, let $V_i$ be a neighborhood of $y(S_i)$ such that $\closure{V_i}\subseteq U_i$.  Let $U$ be a neighborhood of $\pi(y)$ such that $\closure{U}\subseteq\bigcap_{i<n}S_i$.  Set $V=\inv{\pi}U\cap\bigcap_{i<n}\inv{\pi_{S_i}}V_i$.  Then $V$ is a neighborhood of $y$ and we have
\begin{equation*}
\closure{V}\subseteq\bigcap_{i<n}\inv{\pi}S_i\cap\bigcap_{i<n}\inv{\pi_{S_i}}U_i=\bigcap_{i<n}\inv{\pi_{S_i}}U_i;
\end{equation*}
whence, $\closure{V}$ is disjoint from $C$.

Now suppose there is a continuous map $f\colon X\rightarrow [0,1]$ such that $f(\pi(y))=1$ and $f``(X\setminus U)=\{0\}$.  For each $i<n$, likewise suppose there is a continuous map $f_i\colon Y_{S_i}\rightarrow [0,1]$ such that $f_i(y(S_i))=1$ and $f``(Y_{S_i}\setminus U_i)=\{0\}$.  Define $g\colon \bigcap_{i<n}\inv{\pi}S_i\rightarrow [0,1]$ by $z\mapsto f(\pi(z))f_0(z(S_0))\cdots f_{n-1}(z(S_{n-1}))$.  Define $h:\inv{\pi}\bigl(X\setminus\closure{U}\bigr)\rightarrow [0,1]$ by $z\mapsto 0$.  By the pasting lemma, $g\cup h$ is continuous and separates $y$ and $C$.

For (\ref{hereddisconnamalgamative}), suppose $C$ is a nonempty connected subset of $Y$ and $X$ and $Y_S$ are hereditarily disconnected for all $S\in\msS$.  Then $\pi``C$ is connected; whence, $\pi``C=\{p\}$ for some $p\in X$.  For each $S\in\msS$, if $p\in S$, then $\pi_S``C$ is connected; whence, $\card{\pi_S``C}=1$.  Thus, $\card{C}=1$.

For (\ref{zerodimamalgamative}), suppose $S\in\msS$ and $U$ open in $Y_S$ and $y\in\inv{\pi_S}U$.  Let $V$ be a clopen neighborhood of $y(S)$ contained in $U$.  Then $\inv{\pi_S}V$ is clopen in $\inv{\pi}S$.  Let $W$ be a clopen neighborhood of $\pi(y)$ contained in $S$.  Then $\inv{\pi}W\cap\inv{\pi_S}V$ is a clopen neighborhood of $y$ contained in $\inv{\pi_S}U$.

For the converse, first note that each of the classes (\ref{T0amalgamative})\nbd-(\ref{zerodimamalgamative}) is closed with respect to subspaces.  Second, $Y_S$ can be embedded in $Y$ for all $S\in\msS$ because $\prod_{p\in S\in\msS} Y_S$ is a subspace of $Y$ for all $p\in X$.  Finally, $X$ can be embedded in $Y$ because the amalgam of $\la\{f(S)\}\ra_{S\in\msS}$ is homeomorhpic to $X$ for all $f\in\prod_{S\in\msS}Y_S$.
\end{proof}

A countable product of metrizable spaces is metrizable; the next theorem is the analog for amalgams.

\begin{theorem}\label{THMmetricamalgam}
Suppose $X$ and $Y_S$ are metrizable for all $S\in\msS$ and there is a countable $\msT\subseteq\msS$ such that $\card{Y_S}=1$ for all $S\in\msS\setminus\msT$.  Then $Y$ is metrizable.
\end{theorem}
\begin{proof}
Since $Y$ is $T_3$ by Theorem~\ref{THMamalgamative}, it suffices to exhibit a $\sigma$\nbd-locally finite base for $Y$.  For each $T\in\msT$, let $\bigcup_{n<\omega}\mcU_{T,n}$ be a $\sigma$\nbd-locally finite base for $Y_T$; let $\bigcup_{n<\omega}\mcU_n$ be a $\sigma$\nbd-locally finite base for $X$.  For each $n<\omega$ and $\tau\in\Fn{\msT}{\omega}$, set $\mcU_{n,\tau}=\bigl\{U\in\mcU_n: \closure{U}\subseteq\bigcap\dom\tau\bigr\}$ and
\begin{equation*}
\mcV_{n,\tau}=\Bigl\{\inv{\pi}U\cap\bigcap_{T\in\dom\tau}\inv{\pi_T}U_T: U\in\mcU_{n,\tau}\text{ and }(\forall T\in\dom\tau)(U_T\in\mcU_{T,\tau(T)})\Bigr\}.
\end{equation*}
Then $\bigcup_{n<\omega}\bigcup_{\tau\in\Fn{\msT}{\omega}}\mcV_{n,\tau}$ is easily verified to be a $\sigma$\nbd-locally finite base for $Y$.
\end{proof}

In general, productiveness is logically incomparable to amalgamativeness: the class of finite $T_0$ spaces is amalgamative but only finitely productive; the class of powers of $2$ is productive but not amalgamative.  However, all amalgamative classes are finitely productive because if $X\in\msS$ and $\card{Y_S}=1$ for all $S\in\msS\setminus\{X\}$, then $Y\homeo X\times Y_X$.

Given Theorem~\ref{THMamalgamative}, it is tempting to conjecture that amalgams are really subspaces of products in disguise.  This conjecture is false.  To see this, consider the class of nonempty Urysohn spaces.  This class is closed with respect to arbitrary products and subspaces, yet, as demonstrated by the following example, this class is not amalgamative.

\begin{example}\label{EXurysohnNotAmalgamative}
Let $X=\rationals$ with the topology generated by $\{\rationals\setminus K\}$ and the order topology of $\rationals$ where $K=\{2^{-n}: n<\omega\}$.  Then $X$ is Urysohn.  Let $\rationals\setminus K\in\msS$ and, for all $S\in\msS$, let $\card{Y_S}=1$ if $S\not=\rationals\setminus K$.  Set $Y_{\rationals\setminus K}=2$ (with the discrete topology).  Then all the factors of $Y$ are Urysohn.  For each $i<2$, define $y_i\in Y$ by $\{y_i\}=\inv{\pi}\{0\}\cap\inv{\pi_{\rationals\setminus K}}\{i\}$.  Suppose $U_0$ and $U_1$ are disjoint closed neighborhoods of $y_0$ and $y_1$, respectively.  Then $\pi``U_0$ and $\pi``U_1$ are neighborhoods of $0$.  Therefore, $2^{-n}\in\closure{\pi``U_0}\cap\closure{\pi``U_1}$ for some $n<\omega$.  If $2^{-n}\in S\in\msS$, then $\card{Y_S}=1$; hence, $\{\inv{\pi}S: 2^{-n}\in S\in\msS\}$ is a local subbase for $y_2$ where $\{y_2\}=\inv{\pi}\{2^{-n}\}$.  Since $2^{-n}\in\closure{\pi``U_0}\cap\closure{\pi``U_1}$, every finite intersection of elements of this local subbase will intersect $U_0$ and $U_1$.  Hence, $y_2\in\closure{U_0}\cap\closure{U_1}=U_0\cap U_1$, which is absurd.  Therefore, $Y$ is not Urysohn.
\end{example}

In the above example, the base space and all the factors of $Y$ are totally disconnected.  Therefore, no amalgamative class both contains all the nonempty totally disconnected spaces and is contained in the class of nonempty Urysohn spaces.

\begin{question}
Is the class of nonempty realcompact spaces amalgamative?
\end{question}

Despite Example~\ref{EXurysohnNotAmalgamative}, there is a sense in which $Y$ is almost homeomorphic to a subspace of the product of its factors.  For each $S\in\msS$, let $Z_S$ be $Y_S$ with an added point $q_S$ whose only neighborhood is $Z_S$.  Then $Y$ is easily seen to be homeomorphic to the set
\begin{equation*}
\bigcup_{p\in X}\Bigl\{z\in\prod_{S\in\msS}Z_S: (\forall S\in\msS) (z(S)=q_S \Leftrightarrow p\not\in S)\Bigr\}
\end{equation*}
with the subspace topology inherited from $\prod_{S\in\msS}Z_S$.  Moreover, this result still holds if we make $q_S$ isolated for all clopen $S\in\msS$.

Let us make some auxillary definitions relating amalgams to continuous maps and subspaces.  

\begin{definition}
Suppose, for each $S\in\msS$, that $Z_S$ is a nonempty space and $f_S\colon Y_S\rightarrow Z_S$.  Let $Z$ be the amalgam of $\la Z_S\ra_{S\in\msS}$.  Then the \emph{amalgam} of $\la f_S\ra_{S\in\msS}$ is the map $f$ defined by
\begin{equation*}
f=\bigcup_{p\in X}\prod_{p\in S\in\msS} f_S.
\end{equation*}
\end{definition}

In the above definition, it is immediate that $f$ is a map from $Y$ to $Z$.  Moreover, if $f_S$ is continuous for each $S\in\msS$, then $f$ is a continuous map from $Y$ to $Z$.  Similarly, an amalgam of homeomorphisms is a homeomorphism.

\begin{definition}
Suppose $W$ is a subspace of $X$.  The \emph{reduced amalgam} of $\la Y_S\ra_{S\in\msS}$ over $W$ is the space $Z$ defined as follows.  Set $\msT=\{S\cap W: S\in\msS\}\setminus\{\emptyset\}$.  Then $\msT\in\subbase{W}$.  Given $S_0,S_1\in\msS$, declare $S_0\sim S_1$ if $S_0\cap W=S_1\cap W$.  For each $T\in\msT$, let $\varepsilon(T)$ be the unique $\msE$ that is an equivalence class of $\sim$ for which $W\cap\bigcap\msE=T$.  For all $T\in\msT$, set $Z_T=\prod_{S\in\varepsilon(T)}Y_S$.  Let $Z$ be the amalgam of $\la Z_T\ra_{T\in\msT}$.
\end{definition}

In the above definition, $Z$ is homeomorphic to $\bigcup_{p\in W}\prod_{p\in S\in\msS}Y_S$ with the subspace topology inherited from $Y$.

\section{Connectifiable amalgams}

Theorems~\ref{THMamalgamative} and \ref{THMmetricamalgam} demonstrate similarities between products and amalgams.  Of course, amalgams would not be very interesting if there were no major differences between them and products.  Such differences arise for connectedness: unlike a product, an amalgam can be connected even if all its factors are not; connectedness of the base space is sufficient in most cases.  Path\nbd-connectedness of an amalgam with a path\nbd-connected base space is harder to guarantee, but not by much.  Some new positive cennectification results fall out as corollaries.

\begin{theorem}\label{THMconnamalgam}
Suppose $X$ is connected (path\nbd-connected) and there is a finite $E\subseteq X$ such that for all $S\in\msS$ we have $E\not\subseteq S$ or $Y_S$ is connected (path\nbd-connected).  Then $Y$ is connected (path\nbd-connected).
\end{theorem}
\begin{proof}
Let $y_0,y_1\in Y$.  It suffices to show that $y_0$ and $y_1$ are in the same component (path component).  Let $E=\{p_0,\ldots,p_{n-1}\}$ and set $p_n=\pi(y_1)$.  Recursively define $z_0,\ldots,z_{n+1}\in Y$ and $f_0,\ldots,f_n\in\prod_{S\in\msS}Y_S$ as follows.  Set $z_0=y_0$.  Given $z_i$ for some $i<n+1$, choose $f_i$ such that it extends $z_i \cup (y_1\restrict(\dom y_1\setminus\dom z_i))$.  Set $z_{i+1}=f_i\restrict\{S\in\msS:p_i\in S\}$.  For each $i<n+1$, let $C_i$ be the amalgam of $\la\{f_i(S)\}\ra_{S\in\msS}$, which is homeomorphic to $X$ and therefore connected (path\nbd-connected).  Then $z_i, z_{i+1}\in C_i$ for all $i<n+1$; hence, $z_0,\ldots,z_{n+1}$ are in the same component (path component).  

Therefore, it suffices to show that $z_{n+1}$ and $y_1$ are in the same component (path component).  Since $\pi(z_{n+1})=p_n=\pi(y_1)$, we have $\dom z_{n+1}=\dom y_1$. Set $\msA=\dom y_1$ and $\msB=\{S\in\msA: y_1(S)\not=z_{n+1}(S)\}$.  Then $z_{n+1},y_1\in\prod_{S\in\msB}Y_S\times\prod_{S\in\msA\setminus\msB}\{y_1(S)\}$; hence, it suffices to show that $Y_S$ is connected (path\nbd-connected) for all $S\in\msB$.  Suppose $S\in\msB$ and $E\not\subseteq S$.  Then choose the least $i<n$ such that $p_i\not\in S$.  Then $S\not\in\dom z_{i+1}$.  Choose the least $j<n+2$ such that $i+1<j$ and $S\in\dom z_j$.  Then $z_j(S)=f_{j-1}(S)=y_1(S)$ because $S\not\in\dom z_{j-1}$.  Hence, $z_j(S)=y_1(S)\not=z_{n+1}(S)$; hence $j<n+1$.  Choose the least $k<n+2$ such that $j<k$ and $S\in\dom z_k$ and $z_k(S)\not=y_1(S)$.  If $S\not\in\dom z_{k-1}$, then $z_k(S)=y_1(S)$, which is absurd.  Hence, $S\in\dom z_{k-1}$; hence, $z_k(S)=z_{k-1}(S)=y_1(S)$, which is also absurd.  Therefore, for all $S\in\msB$, we have $E\subseteq S$; whence, $Y_S$ is connected (path\nbd-connected).
\end{proof}

\begin{example}\label{EXcone}  Suppose $X=[0,1]$ and $\msS=\{U\subseteq [0,1]: U \text{ open}\}$ and $\card{Y_S}=1$ for all $S\in\msS\setminus\{[0,1)\}$.  Then $Y$ is homeomorphic to the cone over $Y_{[0,1)}$. If $1\in S\in\msS$, then $\card{Y_S}=1$; hence, Theorem~\ref{THMconnamalgam} implies $Y$ is path\nbd-connected.  Thus, Theorem~\ref{THMconnamalgam} may be interpreted as constructing a class of generalized cones.
\end{example}

\begin{corollary}\label{CORpropconncation}
Suppose $i\in\{1,2,3,3\frac{1}{2}\}$ and $X$ has a proper $T_i$ connectification $\tilde{X}$ and $Y_S$ is $T_i$ for all $S\in\msS$.  Then $Y$ has a proper $T_i$ connectification $\tilde{Y}$.  If Moreover, if $\tilde{X}$ is path\nbd-connected, then we may choose $\tilde{Y}$ to be path\nbd-connected.
\end{corollary}
\begin{proof}
Fix $p\in\tilde{X}\setminus X$.  For each $S\in\msS$, let $\Phi(S)$ be an open subset of $\tilde{X}\setminus\{p\}$ such that $\Phi(S)\cap X=S$.  Extend $\Phi``\msS$ to some $\tilde{\msS}\in\subbase{\tilde{X}}$.  For all $S\in\msS$, set $\tilde{Y}_{\Phi(S)}=Y_S$.  For all $S\in\tilde{\msS}\setminus\Phi``\msS$, set $\tilde{Y}_S=1$.  Let $\tilde{Y}$ be the amalgam of $\la\tilde{Y}_S\ra_{S\in\tilde{\msS}}$.  By Theorem~\ref{THMamalgamative}, $\tilde{Y}$ is $T_i$; by Theorem~\ref{THMconnamalgam}, $\tilde{Y}$ is connected, for $\card{\tilde{Y}_S}=1$ if $p\in S\in\tilde{\msS}$.  Define $f:Y\rightarrow\tilde{Y}$ as follows.  Given $y\in Y$, let $\pi(f(y))=\pi(y)$; set $f(y)(\Phi(S))=y(S)$ for all $S\in\dom y$; set $f(y)(S)=0$ for all $S\in\tilde{\msS}\setminus\Phi``\dom y$ such that $\pi(y)\in S$.  Then $f$ is an embedding of $Y$ into $\tilde{Y}$ with dense range $\inv{\pi}X$; hence, $\tilde{Y}$ is a proper $T_i$ connectification of $Y$.  Finally, by Theorem~\ref{THMconnamalgam}, $\tilde{Y}$ is path\nbd-connected if $\tilde{X}$ is.
\end{proof}

The previously known result most similar to Corollary~\ref{CORpropconncation} is due to Druzhinina and Wilson: \cite{druzhinina}: if all the path components of a $T_2$ ($T_3$, metric) space are open and have proper pathwise connectifications, then the space has a $T_2$ ($T_3$, metric) proper pathwise connectification.

\begin{proof}[Proof of Theorem~\ref{THMinfprodinfunion}]
Every infinite product is an infinite product of countably infinite subproducts; every infinite topological sum is a countably infinite topological sum of topological sums.  Moreover, products preserve the property of having a $T_i$ pathwise connectification; topological sums preserve the $T_i$ axiom and metrizability.  Therefore, we only need to prove the theorem for all countably infinite products of countably infinite topological sums.  Set $X=\omega^\omega$ with the product topology.  For each $m,n<\omega$, let $Z_{m,n}$ be a nonempty $T_i$ space and let $S_{m,n}=\{p\in X: p(m)=n\}$; set $Y_{S_{m,n}}=Z_{m,n}$.  Set $\msS=\{S_{m,n}:m,n<\omega\}\in\subbase{X}$.  Then clearly $Y\homeo\prod_{m<\omega}\bigoplus_{n<\omega}Z_{m,n}$.  Since $X\homeo\reals\setminus\rationals$, there is a proper metrizable pathwise connectification of $X$, namely a copy of $\reals$.  By Corollary~\ref{CORpropconncation}, $Y$ has a proper $T_i$ pathwise connectification.  For the metrizable case, construct a connectification $\tilde{Y}$ of $Y$ as in the proof of Corollary~\ref{CORpropconncation}, with $\tilde{X}$ chosen to be homeomorphic to $\reals$.  Since $\msS$ is countable, the space $\tilde{Y}$ is metrizable by Theorem~\ref{THMmetricamalgam}.
\end{proof}

If we care about connectedness but not path\nbd-connectedness, Theorem~\ref{THMconnamalgam} and Corollary~\ref{CORpropconncation} can be considerably strengthened.

\begin{theorem}\label{THMconnamalgam2}
Suppose $X$ is connected and either $X\not\in\msS$ or $Y_X$ is connected.  Then $Y$ is connected.
\end{theorem}
\begin{proof}
Let $y_0, y_1\in Y$.  It suffices to show $y_1$ is in the closure of the component of $y_0$.  Let $U$ be a basic open neighborhood of $y_1$.  Then there exist $n<\omega$ and $\la S_i\ra_{i<n}\in(\dom y_1)^n$ and $\la U_i\ra_{i<n}$ such that $U_i$ is an open neighborhood of $y_1(S_i)$ for all $i<n$ and $U=\bigcap_{i<n}\inv{\pi_{S_i}}U_i$.  Choose $f\in\prod_{S\in\msS}Y_S$ such that $f$ extends $y_0$.  Then there exists $E\subseteq X$ such that $E$ is finite and $E\not\subseteq S$ for all $S\in\{S_i:i<n\}\setminus\{X\}$.  For each $S\in\msS$, set $Z_S=Y_S$ if $Y_S$ is connected or $S\in\{S_i:i<n\}$; otherwise, set $Z_S=\{f(S)\}$.  Let $Z$ be the amalgam of $\la Z_S\ra_{S\in\msS}$.  Then $Z$ is connected by Theorem~\ref{THMconnamalgam}.  Moreover, $y_0\in Z$ and $Z\cap U\not=\emptyset$.  Thus, $y_1$ is in the closure of the component of $y_0$.
\end{proof}

\begin{corollary}\label{CORconncation}
Suppose $i\in\{1,2,3,3\frac{1}{2}\}$ and $X$ has a $T_i$ connectification and $Y_S$ is $T_i$ for all $S\in\msS$.  Further suppose $X$ has a proper $T_i$ connectification or $X\not\in\msS$ or $Y_X$ is connected.  Then $Y$ has a $T_i$ connectification.
\end{corollary}
\begin{proof}
If $X$ has a proper $T_i$ connectification, then so does $Y$ by Corollary~\ref{CORpropconncation}.  If $X$ is $T_i$ and connected but has no proper $T_i$ connectification, then $Y$ is connected by Theorem~\ref{THMconnamalgam2}.
\end{proof}

\section{A large path\nbd-connected homogeneous compactum}

\begin{definition}
We say that a homogeneous compactum is \emph{exceptional} if it is not homeomorphic to a product of dyadic compacta and first countable compacta.
\end{definition}

In the previous section, we constructed a machine for strengthening connectification results.  Next, we construct a machine that takes a homogoeneous compactum and produces a path\nbd-connected homogeneous compactum.  Applying this machine to a particular homogeneous compactum with cellularity $\cardc$, we get a path\nbd-connected homogeneous compactum with cellularity $\cardc$.  Moreover, more careful analysis of the latter space's connectedness properties shows that it is exceptional.

All compact groups are dyadic, and most other known examples of homogeneous compacta are products of first countable compacta (see \cite{Kunen, vanMill}).  Besides the exceptional homogeneous compactum we shall construct, there is, to the best of the author's knowledge, only one known example of an exceptional homogeneous compactum, and its existence is independent of ZFC.  In \cite{vanMill}, van Mill constructed a compactum $K$ satisfying $\piweight{K}=\omega$ and $\character{K}=\omega_1$.  Clearly, $\character{Z}=\omega\leq\piweight{Z}$ for all first countable spaces $Z$.  Moreover, Efimov~\cite{efimov} and Gerlits~\cite{gerlits} independently proved that $\pi\character{Z}=\weight{Z}$ for all dyadic compacta $Z$.  Hence, $\character{Z}\leq\piweight{Z}$ for all $Z$ homeomorphic to products of dyadic compacta and first countable compacta; hence, $K$ is not homeomorphic to such a product.  Under the assumption $\cardp>\omega_1$ (which follows from $\mathrm{MA}+\neg\mathrm{CH}$), van Mill proved that $K$ is homogeneous.  However, van Mill also noted that all homogeneous compacta $Z$ satisfy $2^{\character{Z}}\leq 2^{\piweight{Z}}$ as a corollary of a result of van Douwen \cite{vanDouwen}.  In particular, if $2^\omega<2^{\omega_1}$, then $K$ is not homogeneous.

\begin{definition}
Given a topological space $Z$, let $\Aut{Z}$ denote the group of autohomeomorphisms of $Z$.  Let $\Aut{Z}$ act on $Z$ in the natural way: $gz=g(z)$ for all $z\in Z$ and $g\in\Aut{Z}$.  Let $\Aut{Z}$ act on $\powset{Z}$ such that $gE=g``E$ for all $E\subseteq Z$ and $g\in\Aut{Z}$.
\end{definition}

\begin{lemma}\label{LEMtranstabhomog}
Let $G$ be the stabilizer of $\msS$ in $\Aut{X}$.  Suppose $Z$ is a homogeneous space and $Y_S=Z$ for all $S\in\msS$.  Further suppose $G$ acts transitively on $X$.  Then $Y$ is homogeneous.
\end{lemma}
\begin{proof}
Let $y_0,y_1\in Y$.  Choose $g\in G$ such that $g(\pi(y_0))=\pi(y_1)$.  Define $f\colon Y\rightarrow Y$ as follows.  Given $y\in Y$, let $\pi(f(y))=g(\pi(y))$ and $f(y)(gS)=y(S)$ for all $S\in\dom y$.  Then $f\in\Aut{Y}$ because $f``(\inv{\pi_S}U)=\inv{\pi_{gS}}U$ and $\inv{f}(\inv{\pi_S}U)=\inv{\pi_{\inv{g}S}}U$ for all $S\in\msS$ and $U$ open in $Z$.  Since $y_1,f(y_0)\in Z^{\dom y_1}$, there exists $\la h_S\ra_{S\in\msS}\in\Aut{Z}^{\msS}$ such that $\bigl(\prod_{S\in\dom y_1}h_S\bigr)(f(y_0))=y_1$.  Let $h$ be the amalgam of $\la h_S\ra_{S\in\msS}$.  Then $h\in\Aut{Y}$ and $h(f(y_0))=y_1$.  Thus, $Y$ is homogeneous.
\end{proof}

\begin{lemma}\label{LEMindamalgam}
Suppose $X$ and $Y_S$ are $T_3$ and $\ind Y_S=0$ for all $S\in\msS$.  Then $\ind Y=\ind X$.
\end{lemma}
\begin{proof}
Set $n=\ind X$.  By (\ref{zerodimamalgamative}) of Theorem~\ref{THMamalgamative}, we may assume $n>0$.  We may also assume the lemma holds if $X$ is replaced by a $T_3$ space with small inductive dimension less than $n$.  First, $Y$ is $T_3$ by Theorem~\ref{THMamalgamative}.  Next, given any $f\in\prod_{S\in\msS}Y_S$, the amalgam of $\la\{f(S)\}\ra_{S\in\msS}$ is homeomorphic to $X$; hence, $\ind Y\geq n$.  Let $y\in Y$ and let $U$ be a neighborhood of $y$.  Then $y\in V_0\subseteq U$ where $V_0=\bigcap_{i<m}\inv{\pi_{S_i}}U_i$ for some $m<\omega$ and $\la S_i\ra_{i<m}\in(\dom y)^m$ and $\la U_i\ra_{i<m}$ such that $U_i$ is a clopen neighborhood of $y(S_i)$ for all $i<m$.  Let $W$ be a neighborhood of $\pi(y)$ such that $W\subseteq\bigcap_{i<m}S_i$ and $\ind\bd W<n$.  Set $V_1=V_0\cap\inv{\pi}W$.  

It suffices to show that $\ind\bd V_1<n$.  Set $V_2=\inv{\pi}\bd W$.  Then $\bd V_1=V_0\cap V_2$; hence, it suffices to show that $\ind V_2<n$.  Let $Z$ be the reduced amalgam of $\la Y_S\ra_{S\in\msS}$ over $\bd W$.  Then $Z\homeo V_2$ and $\ind Z=\ind\bd W$ because $\ind\bd W<n$ and every factor of $Z$, being a product of factors of $Y$, is zero\nbd-dimensional.
\end{proof}

\begin{theorem}\label{THMpathconncellc}
There is a path\nbd-connected homogeneous compact Hausdorff space $Y$ with cellularity $\cardc$, weight $\cardc$, and small inductive dimension $1$.  Moreover, $Y$ is not homeomorphic to a product of compacta that all have character less than $\cardc$ or have $\cardc$ a caliber.  In particular, $Y$ is exceptional.
\end{theorem}
\begin{proof}
Let $X$ be the unit circle $\{\la x,y\ra\in\reals^2:x^2+y^2=1\}$.  Let $\msS$ be the set of open semicircles contained in $X$.  Let $\gamma$ be an indecomposable ordinal (\emph{i.e.}, not a sum of two lesser ordinals)  strictly between $\omega$ and $\omega_1$.  For each $S\in\msS$, let $Y_S$ be $2^\gamma$ with the topology induced by its lexicographic ordering.  It is easily seen that $Y_S$ is zero\nbd-dimensional compact Hausdorff and $\weight{Y_S}=\cell{Y_S}=\cardc$.  Moreover, $Y_S$ is homogeneous \cite{Maurice}.  Since $\card{\msS}=\cardc$, we have $\weight{Y}=\cell{Y}=\cardc$.  Moreover, $Y$ is compact Hausdorff by Theorem~\ref{THMamalgamative}.  Since no $S\in\msS$ contains a pair of antipodes, $Y$ is path\nbd-connected by Theorem~\ref{THMconnamalgam}.  The stabilizer of $\msS$ in $\Aut{X}$ contains all the rotations of $X$ and therefore acts transitively on $X$; hence, $Y$ is homogeneous by Lemma~\ref{LEMtranstabhomog}.  Also, by Lemma~\ref{LEMindamalgam}, $\ind Y=\ind X=1$.  

Seeking a contradiction, suppose $Y$ is homeomorphic to a product of compacta that all have character less than $\cardc$ or have $\cardc$ a caliber.  Then there exist a compactum $Z$ with $\cardc$ a caliber, a sequence of nonsingleton compacta $\la W_i\ra_{i\in I}$ all with character less than $\cardc$, and a homeomorphism $\varphi$ from $Z\times\prod_{i\in I}W_i$ to $Y$.  Clearly, $W_i$ is path\nbd-connected for all $i\in I$.  Choose $p\in X$.  Then $\inv{\varphi}\inv{\pi}\{p\}$ is a $G_\delta$\nbd-set; hence, there exist a nonempty $Z_0\subseteq Z$ and $J\in[I]^{\leq\omega}$ and $q\in\prod_{j\in J}W_j$ such that $Z_0\times\{q\}\times\prod_{i\in I\setminus J}W_i\subseteq\inv{\varphi}\inv{\pi}\{p\}$.  Since $\inv{\pi}\{p\}=\prod_{p\in S\in\msS}Y_S$, which is zero\nbd-dimensional, $Z_0\times\{q\}\times\prod_{i\in I\setminus J}W_i$ is also zero\nbd-dimensional; hence, $\prod_{i\in I\setminus J}W_i$ is also zero\nbd-dimensional.  Hence, $W_i$ is not connected for all $i\in I\setminus J$; hence, $I=J$; hence, $I$ is countable.  Set $W=\prod_{j\in J}W_j$.  Then $\character{W}<\cardc$ because $\cf\cardc>\omega$.

Let $H\subseteq X$ be an open arc subtending $\pi/2$ radians.  Set $\msT=\{S\in\msS:H\subseteq S\}$.  Then $\card{\msT}=\cardc$.  Choose a nonempty open box $U\times V\subseteq Z\times W$ such that $U\times V\subseteq\inv{\varphi}\inv{\pi}H$ and $U=\bigcup_{n<\omega}U_n$ where $U_n$ is open and $\closure{U}_n\subseteq U_{n+1}$ for all $n<\omega$.  Choose $r\in V$ and set $\kappa=\character{r,W}<\cardc$.  Let $\la V_\alpha\ra_{\alpha<\kappa}$ enumerate a local base at $r$.  By compactness, we may choose, for each $\alpha<\kappa$ and $n<\omega$, a finite set $\sigma_{n,\alpha}$ of basic open subsets of $Y$ such that $\closure{U}_n\times\{r\}\subseteq\inv{\varphi}\bigcup\sigma_{n,\alpha}\subseteq U_{n+1}\times V_\alpha$.  Set $G=\bigcup_{n<\omega}\bigcap_{\alpha<\kappa}\bigcup\sigma_{n,\alpha}$.  Since $\kappa<\cardc$, there exist nonempty $\msR\subseteq\msT$ and $E\subseteq\bigcup_{x\in H}\prod_{x\in S\in\msS\setminus\msR}Y_S$ such that $G=E\times\prod_{S\in\msR}Y_S$.  Hence, $\cell{G}=\cardc$.  Since $\inv{\varphi}G=U\times\{r\}$, we have $\cell{U}=\cardc$.  Since $U$ is an open subset of $Z$, we have $\cell{Z}\geq\cardc$, which yields our desired contradiction, for $\cardc\in\caliber{Z}$. 
\end{proof}

\begin{remark}
If there is a homogenous compactum with cellularity $\kappa>\cardc$ (that is, if van Douwen's Problem (see \cite{Kunen}) has a positive solution), then the proof of Theorem~\ref{THMpathconncellc} is easily modified to produce a path\nbd-connected homogeneous compactum with cellularity $\kappa$.
\end{remark}

\end{document}